\documentclass[12pt]{amsart}

\usepackage{enumerate}

\numberwithin{equation}{section}

\newcommand{\goto}{\rightarrow}

\newcommand{\bN}{\mathbb{N}}
\newcommand{\bP}{\mathbb{P}}

\newcommand{\bR}{\mathbb{R}}

\newcommand{\bZ}{\mathbb{Z}}

\newcommand{\cF}{\mathcal{F}}

\newcommand{\cP}{\mathcal{P}}

\newcommand{\XX}{\mathbf{X}}
\newcommand{\YY}{\mathbf{Y}}

\newcommand{\xx}{\mathbf{x}}
\newcommand{\yy}{\mathbf{y}}

\newcommand{\bg}{\bar {g}}

\newcommand{\CZZ}{\hat {\mathbb{Z}}}

\theoremstyle{plain} \newtheorem{Theo}{Theorem}[section]
\theoremstyle{plain} 
\theoremstyle{plain} \newtheorem{Cor}[Theo]{Corollary}
\theoremstyle{plain} \newtheorem{Prop}[Theo]{Proposition}
\theoremstyle{remark} \newtheorem{Rem}[Theo]{Remark}
\theoremstyle{definition} 
\theoremstyle{definition} 
\theoremstyle{remark}

\newcommand{\ee}{\mathbf{e}}

\begin{document}

\title[Coalescing Brownian motions]{Balls--in--boxes duality for coalescing random walks
and coalescing Brownian motions}

\author{Steven N. Evans}

\email{evans@stat.Berkeley.EDU}

\address{Department of Statistics \#3860 \\
  University of California at Berkeley \\
367 Evans Hall \\
Berkeley, CA 94720-3860 \\
U.S.A}

\thanks{SNE supported in part by NSF grant DMS-0071468
and a Miller Institute for Basic Research in Science research
professorship, XZ supported by an NSERC grant}

\author{Xiaowen Zhou}

\email{xzhou@mathstat.concordia.ca}

\address{Department of Mathematics and Statistics \\
Concordia University \\
7141 Sherbrooke Street West \\
Montreal, Quebec  \\
H4B 1R6 \\
CANADA}

\begin{abstract}
We present a duality relation between two
systems of coalescing random
walks and an analogous duality relation between two
systems of coalescing
Brownian motions.  Our results extends previous work in
the literature and we apply it
to the study of a system of coalescing Brownian motions with Poisson
immigration.
\end{abstract}

\maketitle

\section{Introduction}

Consider a system of $m$ indexed particles with locations in $\bR$
that evolve as follows.
Each particle moves according to an independent
standard Brownian motion on $\bR$ until two
particles are at the same location. At this moment a {\em
coalescence} occurs and the particle of higher index starts to move
together with the particle of lower index. We say the particle
with higher index is {\em attached} to the particle with lower
index, which is still {\em free}. The particle system then
continues its evolution in the same fashion.  Note that indices
are not essential here, the collection of locations of the particles
is Markovian in its own right, but it will be convenient to think
of the process as taking values in $\bR^m$ rather than
subsets of $\bR$ with at most $m$ elements.  For definiteness,
we will further assume that the particles are indexed in increasing
order of their initial positions: it it clear that the dynamics
preserve this ordering.  Call the resulting Markov process
$\XX = (X_1, \ldots, X_m)$.

The analogous coalescing simple
random walk has many applications. One successful
example is in voter model, which is particularly well
understood because of a duality relation
with the coalescing random walk (see, for example,
\cite{MR80i:60132, MR2001g:60247}).
Similarly, coalescing Brownian motion plays a
key role in analyzing certain  complex interactive stochastic
systems. For example, in \cite{DEFKZ00}
the coalescing Brownian motion  is dual
to the Brownian stepping-stone model in the sense that it
determines the joint ``moments'' of the latter. This interplay
leads to further results on the Brownian stepping-stone
model in \cite{Zh02}. A ``continuous family'' of coalescing
Brownian motions, usually referred to as the Arratia flow, serves as a
fundamental example in the theory of stochastic flows. See \cite{Ar79, Ha84}
for accounts of this topic.  The Arratia flow is
an example of an interesting noise that is not
generated by Brownian motions or Poisson processes \cite{MR99i:60082, Tsi04}.
More general ``sticky'' flows have recently been considered in
\cite{MR2003d:60114, MR1969575, MR2004c:60180}.

Closed form analytic expressions for features
of the joint distribution of coalescing
Brownian motion are rarely known, but some intriguing
relationships have been observed for stochastic systems involving
coalescing Brownian motions. A self-duality
relation for the Arratia flow is described
in the Introduction of \cite{Ar79}, where the borders
between clusters (that is, pre-images of particles)
are shown to have the same joint
distribution as the locations of particles.
A duality between a system of coalescing Brownain motions
and a system of annihilating Brownian motions is established
in \cite{DEFKZ00}.
A dual relationship is presented in
\cite{STW00} between two system of Brownian motions, in which one
system runs forward in time, the other runs backward in time,
Brownian motions from the same system coalesce and Brownian
motions from different systems reflect on each other. Another
result  along this line is obtained in in \cite{ToWe97}, which
involves a duality on two flows of  Brownian motions moving at
opposite directions of the time interval $(-\infty, \infty)$.
Within each flow, the Brownian motions  coalesce, and meanwhile
each Brownian motion is either reflected or absorbed at $0$
depending on when it reaches $0$. There is no interaction between
the two flows.

The distribution  of $\XX(t)$ is uniquely specified by knowing
for each choice of $y_1 < y_2 < \ldots < y_n$ the joint probabilities
of which ``balls'' $X_1(t), X_2(t), \ldots, X_m(t)$
lie which of the
``boxes'' $[y_1, y_2], [y_2, y_3], \ldots, [y_{n-1}, y_n]$.
That is, the distribution of $\XX(t)$
is determined by the joint distribution of the indicators
\[
I_{ij}^\rightarrow(t,\yy) := 1\{X_i(t) \in [y_{j}, y_{j+1}]\}
\]
for $1 \le i \le m$ and $1 \le j \le n-1$.

Suppose now that $\YY := (Y_1, \ldots, Y_n)$ is another
coalescing Brownian motion.  The distribution of $\YY(t)$
is uniquely specified by knowing for each choice of
$x_1 < x_2 < \ldots < x_n$ the distribution of the indicators
\[
I_{ij}^\leftarrow(t,\xx) := 1\{x_i \in [Y_{j}(t), Y_{j+1}(t)]\}
\]
for $1 \le i \le m$ and $1 \le j \le n-1$.

Thus we can think of a coalescing Brownian motion as being
a set of evolving balls with the distribution at time
$t$ determined by how the balls fall in
a fixed set of boxes, or we
can think of a coalescing Brownian motion as giving
a set of evolving boxes with the distribution at time
$t$ determined by how these boxes contain a fixed set of
balls.  We show that these two points of view are
{\em dual} to each other in the sense that if
$\XX(0) = \xx$ and $\YY(0) = \yy$, then for
each $t \ge 0$ the arrays
of indicators $(I_{ij}^\rightarrow(t,\yy))$
and $(I_{ij}^\leftarrow(t,\xx))$ have the same
joint distribution.  We derive this duality from
an analogous, essentially combinatorial, fact about coalescing simple random walk.

Special cases of the above mentioned duality were proved earlier
in \cite{XiZh03}. Instead of using a discrete approximation approach,
the results there were directly obtained from coalescing Brownian
motions, and, as a result, the proofs were rather lengthy.

Moreover, we extend the Brownian motion result to a situation
where the ``balls'' and the ``boxes'' are allowed to originate at
different points in time (rather than all originating at time
$0$). This latter extended result is then used to analyse the
asymptotic behaviour of a system of coalescing Brownian particles
in which new particles arise according to a homogeneous
space--time Poisson point process.

The rest of this paper is arranged as follows. Section \ref{rw}
contains the preparation and the proof of our main result on the
duality between two coalescing simple random walks. In section
\ref{bm} we generalize this dual relationship to coalescing
Brownian motions starting from possibly different times. Some
known results are re-derived. In section \ref{bmpoi} we further
generalize the duality to a model with both coalescing Brownian
motion and Poisson migration included.

\section{Coalescing random walk}
\label{rw}

A $p$-simple random walk on $\bZ$ is a continuous time simple
random walk that makes jumps at unit rate, and when it makes a
jump from some site it jumps to the right neighbour with
probability $p$ and to the left neighbour with probability $1-p$.
An $m$-dimensional $p$-simple coalescing random walk is defined in
the same way as the coalescing Brownian motion of the
Introduction. When $p=\frac{1}{2}$ we just call this particle
system a simple coalescing random walk.

Some notation is useful to keep track of the interactions among
the particles in the coalescing system.
Let $\cP_m$ denote the set of {\em interval partitions} of the
totality of indices $\bN_m := \{1, \ldots, n\}$. That is, an
element $\pi$ of $\cP_n$ is a collection $\pi = \{A_1(\pi),
\ldots, A_h(\pi)\}$ of disjoint subsets of $\bN_m$ such that
$\bigcup_i A_i(\pi) = \bN_m$ and $a<b$ for all $a\in A_i$, $b\in
A_j$, $i<j$. The sets $A_1(\pi), \ldots A_h(\pi)$ consisting of
consecutive indices are the {\em
intervals} of the partition $\pi$. The integer $h$ is the {\em
length} of $\pi$ and is denoted by $l(\pi)$. Equivalently, we can
think of $\cP_m$ as a set of equivalence relations on $\bN_m$
and write $i \sim_\pi j$ if $i$ and $j$ belong to the same interval
of $\pi \in \cP_m$. Of course, if $i \sim_\pi j$, then
$i \sim_\pi k \sim_\pi j$ for all $i \le k \le j$.

Given $\pi\in \cP_m$, define
$$\alpha_i(\pi) := \min A_i(\pi)$$
to be the left-hand end-point of the $i^{\mathrm{th}}$ interval $A_i(\pi)$.
Put
\[
\bZ^m_\pi := \{(x_1, \ldots, x_m) \in \bZ^m : x_1 \le \ldots \le x_m \text{ and } x_i = x_j \text{ if } i \sim_\pi j\}
\]
and
\[
\begin{split}
\CZZ_\pi^m:=
\{(x_1, \ldots, x_m) \in \bZ^m : x_1 \le \ldots \le x_m \text{ and } x_i = x_j \text{ if and only if } i \sim_\pi j\}.\\
\end{split}
\]
Note that $\bZ^m$ is the
disjoint union of the sets $\CZZ^m_\pi$, $\pi\in\cP_m$.

Write $\XX=(X_1,\ldots,X_m)$ for the coalescing random walk.
If $\XX(t) \in \CZZ^m_\pi$, then the free
particles at time $t$ have indices $\alpha_1(\pi),\ldots,\alpha_{l(\pi)}(\pi)$
and the $i^{\mathrm{th}}$ particle at time $t$ is attached to the free particle with index
\[
\min\{j : 1\le j \le m, \, j \sim_\pi i\}
=
\max\{\alpha_k(\pi) : \alpha_k(\pi) \le i\}.
\]

In order to write down the generator of $\XX$,
we require a final piece of notation.
 Let  $\{\ee_{i}^k: 1\leq i\leq k\}$ be the
set of coordinate vectors in $\bZ^k$; that is, $\ee_{i}^k$ is the
vector that has $i^{\mathrm{th}}$ coordinate $1$ and all other
coordinates  $0$. For $\pi \in \cP_m$, define a map
$K_\pi:\bZ^m_\pi \rightarrow \bZ^{l(\pi)}$ by
\[
K_\pi(\xx)=K_\pi(x_1,\ldots,x_m):=(x_{\alpha_1(\pi)},\ldots,x_{\alpha_{l(\pi)}(\pi)})
\]
Notice that $K_\pi$ is a bijection between $\bZ_\pi^m$ and $\{x
\in \bZ^{ l(\pi)} : x_1 \le x_2 \le \ldots \le x_{l(\pi)}\}$. For
brevity, we will sometimes write $\xx_\pi$ for $K_\pi(\xx)$.

Write $B(\bZ^m)$ for the collection of all bounded functions on
$\bZ^m$. The generator $G$ of $\XX$ is the operator
$G: B(\bZ^m)\rightarrow B(\bZ^m)$ given by
\begin{equation*}
\begin{split}
G f (\xx)
&:=p \sum_{i=1}^{l(\pi)} f\circ
K_\pi^{-1}(\xx_\pi+\ee_{i}^{l(\pi)})
+(1-p) \sum_{i=1}^{l(\pi)}f\circ
K_\pi^{-1}(\xx_\pi-\ee_{i}^{l(\pi)}) \\
& \quad - l(\pi)f\circ
K_\pi^{-1}(\xx_\pi),
\quad f\in B(\bZ^m), \; \xx\in \CZZ_\pi^m, \; \pi\in\cP_m.\\
\end{split}
\end{equation*}
This expression is well-defined, because if
$\xx \in \CZZ_\pi^m$, then
$\xx_\pi$, $\xx_\pi+\ee_{i}^{l(\pi)}$ and $\xx_\pi-\ee_{i}^{l(\pi)}$
are all in
$\{x \in \bZ^{ l(\pi)} : x_1 \le x_2 \le \ldots \le x_{l(\pi)}\}$.

\noindent
{\bf Note:} From now on we will suppress the dependence on dimension
and write $\ee_{i}^{l(\pi)}$ as $\ee_{i}$.

Write $\bZ':=\bZ+\frac{1}{2}=\{i+\frac{1}{2}:i\in \bZ\}$. An
$n$-dimensional $q$-simple coalescing random walk on ${\bZ'}^n$
and its generator $H$
can be defined in the obvious way. Such a
process, with $q=1-p$, will serve as the process dual
to the $p$-simple coalescing random walk on $\bZ^m$
in the following way.

Fix $\xx \in \bZ^m$ with $x_1 \le \ldots \le x_m$
and $\yy \in {\bZ'}^n$ with $y_1 \le \ldots \le y_n$.
By analogy with the notation introduced in the Introduction, put
\[
I_{ij}^\rightarrow(t,\yy) := 1\{X_i(t) \in [y_{j}, y_{j+1}]\}
\]
and
\[
I_{ij}^\leftarrow(t,\xx) := 1\{x_i \in [Y_{j}(t), Y_{j+1}(t)]\}
\]
for $1 \le i \le m$ and $1 \le j \le n-1$.

\begin{Theo}
\label{Tha}
Suppose in the notation above that
$\XX=(X_1,\ldots,X_m)$ is an $m$-dimensional $p$-simple coalescing
random walk starting at $\xx=(x_1,\ldots,x_m)$ with
$x_1\leq \ldots \leq x_m$ and $\YY=(Y_1,\ldots,Y_n)$ is a
$n$-dimensional $(1-p)$-simple coalescing random walk starting at
$\yy=(y_1,\ldots,y_n)$ with $y_1\le \ldots \le y_m$.
Then for each $t \ge 0$
the joint distribution of the $m \times (n-1)$-dimensional
random array $(I_{ij}^\rightarrow(t,\yy))$
coincides with that of the
$m \times (n-1)$-dimensional
random array
$(I_{ij}^\leftarrow(t,\xx))$.
\end{Theo}

\begin{proof}
For a function $g: \{0,1\}^{m(n-1)} \rightarrow \bR$,
a vector $\tilde \xx \in \bZ^m$ with
$\tilde x_1 \le \ldots \le \tilde x_m$,
and a vector
$\tilde \yy \in {\bZ'}^n$ with
$\tilde y_1 \le \ldots \le \tilde y_n$,
set
\begin{equation*}
\bg(\tilde \xx; \tilde \yy)
:=g(1_{[\tilde y_1,\tilde y_2]}(\tilde x_1),
\ldots,
1_{[\tilde y_{n-1},\tilde y_n]}(\tilde x_1),
\ldots,1_{[\tilde y_1, \tilde y_2]}( \tilde x_m),
\ldots,1_{[\tilde y_{n-1}, \tilde y_n]}(\tilde x_m)).
\end{equation*}
We may assume that $\XX$ and $\YY$ are defined
on the same probability space $(\Omega, \cF, \bP)$.
We need to show that
\begin{equation}\label{ta1}
\bP[\bg(\XX_t;\yy)]=\bP[\bg(\xx;\YY_t)].
\end{equation}

For $\tilde \xx\in \bZ^m$, put
$\bg_{\tilde \xx}(\cdot):=\bg(\tilde \xx;\cdot)$,
and for $\tilde \yy'\in \bZ'^n$, put
$\bg_{\tilde \yy}(\cdot):=\bg(\cdot; \tilde \yy)$.
In order to establish (\ref{ta1}),
it suffices by a standard argument
(cf. Section 4.4 in \cite{EtKu86})
to show that
\begin{equation}
\label{genduality}
G(\bg_{\yy})(\xx) = H(\bg_{\xx})(\yy)
\end{equation}
(recall that $G$ and $H$ are the generators of $\XX$ and $\YY$,
respectively).

Fix
$\xx\in\CZZ_{\pi}^m$ and $\yy\in\CZZ_{\varpi}'^n$
for some
$\pi\in\cP_m$ and $\varpi\in\cP_n$.
Set
\[
I^+:=\{i: 1\leq i\leq l(\pi), \;
x_{\alpha_i(\pi)}+\frac{1}{2}
=y_{\alpha_j(\varpi)}
\text{\, for some}\, 1\leq j\leq l(\varpi)\}
\]
and
\[
I^-:=\{i: 1\leq i\leq l(\pi), \; x_{\alpha_i(\pi)}-\frac{1}{2}=y_{\alpha_j(\varpi)} \text{\, for some}\, 1\leq j\leq l(\varpi)\}.\]
Similarly, put
\[
J^-:=\{j: 1\leq j\leq l(\varpi), \; y_{\alpha_j(\varpi)}-\frac{1}{2}=x_{\alpha_i(\pi)} \text{\, for some}\, 1\leq i\leq l(\pi)\}
\]
and
\[J^+:=\{j: 1\leq j\leq l(\varpi), \; y_{\alpha_j(\varpi)}+\frac{1}{2}=x_{\alpha_i(\pi)} \text{\, for some}\, 1\leq i\leq l(\pi)\}.\]

Recall that $x_{\alpha_1(\pi)} < \ldots <
x_{\alpha_{l(\pi)}(\pi)}$ and $y_{\alpha_1(\varpi)} < \ldots <
y_{\alpha_{l(\varpi)}(\varpi)}$. Therefore, for each $i \in I^+$
there is a unique $j \in J^-$ such that
$x_{\alpha_i(\pi)}+\frac{1}{2}=y_{\alpha_j(\varpi)}$ and {\em vice
versa}. Fix such a pair $(i,j)$. Observe that
\[\xx':=\xx+\sum_{k\in A_i(\pi)}\ee_{k}^m=K^{-1}(\xx_\pi+\ee_i)  \]
and
\[\yy':=\yy-\sum_{k\in A_j(\varpi)}\ee_{k}^n=K^{-1}(\yy_\varpi-\ee_j).  \]
Writing $1(B)(\cdot)$ for the indicator function of
a set $B$, we are going to verify that
\begin{equation}\label{ta2}
\left(1([y_{j'}, y_{j'+1}])(x'_{i'})\right)=(1([y'_{j'},
y'_{j'+1}])(x_{i'}))
\end{equation}
by considering all the possible scenarios.

Given any $i'\in A_i(\pi)$ we have:
\begin{itemize}
\item
for $j'=\alpha_j(\varpi)-1$,
\begin{equation*}
\begin{split}
&1([y_{j'},y_{j'+1}])(x'_{i'})=1([y_{j'},y_{j'+1}])(x_{i'}+1) \\
&\quad=0\\
&\quad=1([y_{j'},y_{j'+1}-1])(x_{i'})=1([y'_{j'},y'_{j'+1}])(x_{i'}),
\end{split}
\end{equation*}
\item
for $\alpha_j(\varpi)\leq j'<\max A_j(\varpi)$,
\begin{equation*}
\begin{split}
&1([y_{j'},y_{j'+1}])(x'_{i'})=1([y_{j'},y_{j'+1}])(x_{i'}+1)\\
&\quad=0\\
&\quad=1([y_{j'}-1,y_{j'+1}-1])(x_{i'})
=1([y'_{j'},y'_{j'+1}])(x_{i'}),
\end{split}
\end{equation*}
\item for $j'=\max A_j(\varpi)$,
\begin{equation*}
\begin{split}
&1([y_{j'},y_{j'+1}])(x'_{i'})=1([y_{j'},y_{j'+1}])(x_{i'}+1)\\
&\quad=1\\
&\quad=1([y_{j'}-1,y_{j'+1}])(x_{i'})=1([y'_{j'},y'_{j'+1}])(x_{i'}),
\end{split}
\end{equation*}
\item
and for $j'< \alpha_j(\varpi)-1$ or $j'> \max A_j(\varpi)$,
\begin{equation*}
\begin{split}
&1([y_{j'},y_{j'+1}])(x'_{i'})=1([y_{j'},y_{j'+1}])(x'_{i'})=1([y_{j'},y_{j'+1}])(x_{i'}+1)\\
&\quad=0\\
&\quad=1([y_{j'},y_{j'+1}])(x_{i'})=1([y'_{j'},y'_{j'+1}])(x_{i'}).
\end{split}
\end{equation*}
\end{itemize}
Moreover, given any
$i'\not\in A_i(\pi)$, we have $x_{i'}\neq x_{\alpha_i(\pi)}$.
Hence
\begin{itemize}
\item for $j'=\alpha_j(\varpi)-1$,
\[1([y_{j'},y_{j'+1}])(x'_{i'})=1([y_{j'},y_{j'+1}])(x_{i'})=1([y_{j'},y_{j'+1}-1])(x_{i'})
=1([y'_{j'},y'_{j'+1}])(x_{i'}),\]
\item for $j'=\max A_j(\varpi)$,
\[1([y_{j'},y_{j'+1}])(x'_{i'})=1([y_{j'},y_{j'+1}])(x_{i'})=1([y_{j'}-1,y_{j'+1}])(x_{i'})
=1([y'_{j'},y'_{j'+1}])(x_{i'}),\]
\item for $\alpha_j(\varpi)\leq
 j'< \max A_j(\varpi)$,
\[1([y_{j'},y_{j'+1}])(x'_{i'})=1([y_{j'},y_{j'+1}])(x_{i'})=1([y_{j'}-1,y_{j'+1}-1])(x_{i'})
=1([y'_{j'},y'_{j'+1}])(x_{i'}),\] \item and for
$j'<\alpha_j(\varpi)-1$ or $j'>\max A_j(\varpi)$,
\[1([y_{j'},y_{j'+1}])(x'_{i'})=1([y_{j'},y_{j'+1}])(x_{i'})
=1([y'_{j'},y'_{j'+1}])(x_{i'}).\]
\end{itemize}
Combining the above observations yield (\ref{ta2}).

Therefore,
\[\bg_\yy \circ
K_\pi^{-1}(\xx_{\pi}+\ee_i)=\bg_\xx \circ
K_\pi^{-1}(\yy_{\varpi}-\ee_j).\]
Furthermore, it is easy to see
for $i'\not\in I^+$ that
\[\bg_\yy \circ K_\pi^{-1}(\xx_{\pi}+\ee_{i'})
= \bg_\yy \circ K_\pi^{-1}(\xx_{\pi})\]
and for $j'\not\in J^-$ that
\[\bg_\xx \circ K_\varpi^{-1}(\yy_{\varpi}-\ee_{j'})
=\bg_\xx \circ K_\varpi^{-1}(\yy_{\varpi}).\]

Similarly, for any $i\in I^-$ there exists a unique $j\in
J^+$ such that
$x_{\alpha_i(\pi)}-\frac{1}{2}=y_{\alpha_j(\varpi)}$
and {\em vice versa}.  For such a pair $(i,j)$ we have
\[\bg_\yy \circ K_\pi^{-1}(\xx_{\pi}-\ee_i)=\bg_\xx \circ
K_\pi^{-1}(\yy_{\varpi}+\ee_j).\]
Furthermore, we see for
$i'\not\in I^-$ that
\[\bg_\yy \circ K_\pi^{-1}(\xx_{\pi}-\ee_{i'})= \bg_\yy \circ K_\pi^{-1}(\xx_{\pi}) \]
and for
$j'\not\in J^+$ that
\[\bg_\xx \circ K_\varpi^{-1}(\yy_{\varpi}+\ee_{j'})=\bg_\xx \circ K_\varpi^{-1}(\yy_{\varpi}).\]

Lastly, note that
\[
\bg_\yy \circ K_\pi^{-1}(\xx_{\pi}) = \bg(\xx;\yy)= \bg_\xx \circ
K_\varpi^{-1}(\yy_{\varpi})
\]
and so
\begin{equation*}
\begin{split}
&G(\bg_{\yy})(\xx)-H(\bg_{\xx})(\yy)\\
&\quad=p \sum_{i=1}^{l(\pi)} \left(\bg_{\yy}\circ
K_\pi^{-1}(\xx_\pi+\ee_i)
-\bg_{\yy}\circ K_\pi^{-1}(\xx_\pi)\right)\\
&\qquad + (1-p) \sum_{i=1}^{l(\pi)}\left(\bg_{\yy}\circ
K_\pi^{-1}(\xx_\pi-\ee_i)
-\bg_{\yy}\circ K_\pi^{-1}(\xx_\pi)\right)\\
&\qquad -p \sum_{j=1}^{l(\varpi)} \left(\bg_{\xx}\circ
K_\varpi^{-1}(\yy_{\varpi}-\ee_i)-\bg_{\xx}\circ K_\varpi^{-1}(\yy_\varpi)\right)\\
&\qquad -(1-p) \sum_{j=1}^{l(\varpi)} \left(\bg_{\xx}\circ
K_\varpi^{-1}(\yy_\varpi+\ee_i)-\bg_{\xx}\circ K_\varpi^{-1}(\yy_\varpi)\right)\\
&\quad=p\sum_{i\in I^+}\bg_{\yy}\circ
K_\pi^{-1}(\xx_\pi+\ee_i)
-p\sum_{j\in J^-}\bg_{\xx}\circ K_\varpi^{-1}(\yy_\varpi-\ee_j)\\
&\qquad +(1-p)\sum_{i\in I^-}\bg_{\yy}\circ
K_\pi^{-1}(\xx_\pi-\ee_i)
-(1-p)\sum_{j\in J^+}\bg_{\xx}\circ K_\varpi^{-1}(\yy_\varpi+\ee_j)\\
&\quad=0,\\
\end{split}
\end{equation*}
as required.
\end{proof}

\begin{Rem}
One can see from the proof that it is crucial that the random walks
make only nearest neighbor jumps.
\end{Rem}

\section{Coalescing Brownian motion}
\label{bm}

In this section we will show that the duality in Theorem \ref{Tha}
also holds when the coalescing random walks are replaced by
coalescing Brownian motions. Coalescing Brownian motion can be
defined similarly to  coalescing random walk.
This duality between two coalescing Brownian motions follows if one
can show the unsurprising fact that a coalescing random walk
scaled in time and space in the usual way converges weakly to a
coalescing Brownian motion.

\begin{Prop}
\label{Pro2a} The conclusion of Theorem \ref{Tha} holds when the
coalescing random walks are replaced by coalescing Brownian
motions in the definition of $(I_{ij}^\rightarrow(t,\yy))$ and
$(I_{ij}^\leftarrow(t,\xx))$.
\end{Prop}

We omit the the proof, but remark that
a particularly straightforward martingale argument
proof of the convergence of coalescing
random walk to coalescing Brownian
motion can be given using the following
result that parallels L\'evy's celebrated martingale
characterization of Brownian motion (and is a fairly
simple consequence of that result). We also omit the
proof of this theorem.

\begin{Theo}
\label{Thb} Let $\XX$ be an $m$-dimensional continuous process
with $\XX(0) = \xx$, where $x_1 \le \ldots \le x_m$, and let $\cF^\XX$
denote the filtration generated by $\XX$. Then the following are
equivalent.
\begin{enumerate}[(i)]
\item The process $\XX$ is a coalescing Brownian motion. \item For
each $1\leq i\leq m$, the process $X_i$ is a Brownian motion with
respect to $\cF^\XX$, and for each pair $ 1\leq i<j\leq n$, the
process $\frac{1}{\sqrt 2}(X_j-X_i)$ is a Brownian motion stopped
at $0$ with respect to $\cF^\XX$. \item The process $\XX$ is a
continuous martingale with quadratic variation $\langle
X_i,X_j\rangle_t=t-T_{ij}\wedge t $, where $T_{ij}:=\inf \{s\geq
0:X_i(s)=X_j(s)\}$, $1\leq i \leq j\leq n$.
\end{enumerate}
\end{Theo}

For a function $g: \{0,1\}^{m(n-1)} \rightarrow \bR$,
a vector $\tilde \xx \in \bZ^m$ with
$\tilde x_1 \le \ldots \le \tilde x_m$,
and a vector
$\tilde \yy \in {\bZ'}^n$ with
$\tilde y_1 \le \ldots \le \tilde y_n$,
set
\begin{equation*}
\bg(\tilde \xx; \tilde \yy)
:=g(1_{[\tilde y_1,\tilde y_2]}(\tilde x_1),
\ldots,
1_{[\tilde y_{n-1},\tilde y_n]}(\tilde x_1),
\ldots,1_{[\tilde y_1, \tilde y_2]}( \tilde x_m),
\ldots,1_{[\tilde y_{n-1}, \tilde y_n]}(\tilde x_m)).
\end{equation*}
Proposition \ref{Pro2a} says that for any $\xx$ and $\yy$
\begin{equation*}
\bP[\bg(\XX_t;\yy)]=\bP[\bg(\xx;\YY_t)].
\end{equation*}
By choosing the right function $g$ we can recover some known
dualities. For example, given $\pi=(A_1,\ldots,A_h)\in\cP_n$ and
$y_1<\ldots<y_{2h}$, put
\[\bg(\xx;\yy)=\prod_{j=1}^{h}\prod_{i\in
A_j}1_{[y_{2j-1},y_{2j}]}(x_i), \; \xx\in\bR^n, \, \yy\in\bR^m.\]
 Then Proposition \ref{Pro2a} implies that
\[ \bP\left\{\bigcap_{j=1}^{ h}\bigcap_{i\in A_j}\{X_i(t)\in
[y_{2j-1},y_{2j}]\}\right\}
=\bP\left\{\bigcap_{j=1}^{h}\bigcap_{i\in A_j}\{x_i\in
[Y_{2j-1}(t), Y_{2j}(t)]\} \right\},\]
 which gives the duality in Theorem 1.1 of
\cite{XiZh03}. If we choose
\[\bg(\xx;\yy)=\prod_{i=1}^{n}\left(1-\prod_{j=1}^m
(1-1_{[y_{2j-1},y_{2j}]}(x_i))\right), \; \xx\in\bR^n,
y_1<\ldots<y_{2m},\] then
\[\bP\left\{\bigcap_{i=1}^n\{X_i(t)\in\bigcup_{j=1}^m[y_{2j-1},y_{2j}]\}\right\}=
\bP\left\{\bigcap_{i=1}^n\{X_i(0)\in\bigcup_{j=1}^m[Y_{2j-1}(t),Y_{2j}(t)]\}\right\}.\]
Therefore, Proposition 3.7 in \cite{XiZh03} follows readily.

The duality Proposition \ref{Pro2a} can be generalized to one
involving coalescing Brownian motion starting from different
times.  In order to state this result, it will be convenient to
think of coalescing Brownian motion a little differently from what
we have done so far. As we have defined it, the coalescing
Brownian motion $\XX$ takes values in the space $\{\xx \in \bR^m :
x_1 \le x_2 \le \ldots \le x_m\}$. It will be more convenient to
work with a related process for which we don't impose this
condition. Given an arbitrary $\xx \in \bR^m$, let $\sigma$ be any
permutation of the indices $\{1,2, \ldots, m\}$ such that
$x_{\sigma(1)} \le x_{\sigma(2)} \le \ldots x_{\sigma(m)}$. Let
$\tilde \XX$ be an $\bR^m$-valued process that has the same
distribution as the process $(X_{\sigma^{-1}(1)}, X_{\sigma^{-1}(2)},
\ldots, X_{\sigma^{-1}(m)})$, where $\XX(0) = (x_{\sigma(1)},
x_{\sigma(2)}, \ldots, x_{\sigma(m)})$. It is not difficult to see
that $\tilde \XX$ is a time-homogeneous strong Markov process. The
following result is obvious.

\begin{Cor}
\label{CorPro2a} The duality in Proposition \ref{Pro2a} holds when
the ordered coalescing Brownian motion $\XX$ is replaced by the
unordered coalescing Brownian motion $\tilde \XX$.
\end{Cor}

Given $((s_1, x_1), \ldots, (s_m, x_m)) \in (\bR_+ \times \bR)^m$
with $0 \le s_1 \le s_2 \le \ldots \le s_m$, define a process
$\bar \XX$  taking values in $\{\epsilon\} \cup \bigcup_{k=1}^m
(\bR_+ \times \bR)^k$, where $\epsilon$ is the null vector, as
follows. Let $0 \le \sigma_1 < \ldots < \sigma_\ell$ denote the
distinct elements of $(s_1, \ldots, s_m)$ written in order.  For
$t \in [0, \sigma_1[$, $\bar \XX(t) = \epsilon$. For $t \in
[\sigma_1, \sigma_2[$, $\bar \XX$ evolves as $\tilde \XX(\cdot -
\sigma_1)$ under the initial condition $\tilde \XX(0) = (x_i : s_i
= \sigma_1)$. Inductively, if $\bar \XX(t)$ has been defined on
$[0, \sigma_h[$, then for $t \in [\sigma_h, \sigma_{h+1}[$ (with
the convention $\sigma_{\ell+1} = \infty$), $\bar \XX$ evolves
conditionally independently of $\{\bar \XX(u) : u \in
[0,\sigma_h[\}$ given $\bar \XX(u-)$ as $\tilde \XX(\cdot -
\sigma_h)$ under the initial condition $\tilde \XX(0) = \bar
\XX(u-) \cup (x_i : s_i = \sigma_h)$ (where $\cup$ denotes the
operation of appending one vector to the end of the other).

The following result is a straightforward consequence of
Proposition \ref{Pro2a} and repeated applications of the Markov
property.

\begin{Prop}
\label{Pro2b} Let
$((s_1, x_1), \ldots, (s_m, x_m))$ and $\bar \XX$ be as above, and
let $\YY=(Y_1,\ldots,Y_n)$ be
a coalescing Brownian motion starting at
$\yy=(y_1,\ldots,y_n)$, with $y_1 \le \ldots \le y_m$.
Then, for $t \ge \max_i s_i$, the $m \times (n-1)$-dimensional
random array
\begin{equation*}
\Bigl( 1_{[y_j,y_{j+1}]}(\bar X_i(t)) \Bigr)
\end{equation*}
has the same distribution as
\begin{equation*}
\Bigl( 1_{[Y_j(t-s_i),Y_{j+1}(t-s_i)]}(x_i) \Bigr).
\end{equation*}
\end{Prop}

Given two functions $f,g : \bR_+ \rightarrow \bR$ with
$f(t) \le g(t)$ for all $t$, let $D_t^\rightarrow (f,g)\subset
[0,t]\times\bR$ denote the region sandwiched between the graphs of
$f$ and $g$ up to time $t$.  That is,
 \[
D_t^\rightarrow(f,g):=\{(s,y): 0\leq s\leq t,
 f(s)<y<g(s)\}.
\]
Let $D_t^\leftarrow(f,g):=\{(t-s,y):(s,y)\in D_t(f,g)\}$ be
the region $D_t^\rightarrow (f,g)$ time-reversed at time $t$.
The
conclusion of in Proposition \ref{Pro2b} is that
the random array
\begin{equation*}
\left(1_{[y_j,y_{j+1}]}(\tilde X_i(t))\right)
\end{equation*}
has the same distribution as the random array
\begin{equation*}
\left(1\{(s_i,x_i)\in D_t^\leftarrow (Y_j,Y_{j+1}) \} \right).
\end{equation*}

\section{Coalescing Brownian motion with Poisson migration}
\label{bmpoi}

In this section we are going to study a particle system which can
be described intuitively as follows.  Given a time-space
Poisson random measure $\Pi^+$ on $\bR_+ \times \bR$ with intensity
measure $\lambda \times {\mathrm{Lebesgue}}$.
Particles appear at the atoms of $\Pi^+$.
Once a particle appears, it starts to move. The
existing particles execute
coalescing Brownian motion with possibly different initial times.
Define a set-valued process $S$ by taking $S_t$ to be the
set of locations of those particles at time $t>0$.

The easiest way to define $S$ formally is via the coalescing
Brownian flow $\phi$ of Arratia \cite{Ar79}.  Here $\phi(s,t,x)$
for $s,t,x \in \bR$ with $s \le t$ is a collection of random
variables with the properties
\begin{itemize}
\item the random map $(s,t,x) \mapsto \phi(s,t,x)$ is jointly
measurable, \item for each $s$ and $x$, the map $t \mapsto
\phi(s,t,x)$, $t \ge s$, is continuous, \item for each $s$ and $t$
with $s \le t$, the map $x \mapsto \phi(s,t,x)$ is non-decreasing
and right-continuous, \item for $s \le t \le u$, $\phi(t,u,\cdot)
\circ \phi(s,t,\cdot) = \phi(s,u,\cdot)$, \item for $u \in \bR$,
$(s,t,x) \mapsto \phi(s+u, t+u, x)$ has the same distribution as
$\phi$, \item for $x_1, \ldots, x_m \in \bR$ the process
$(\phi(0,t,x_1), \ldots, \phi(0,t,x_m))_{t \ge 0}$ has the same
distribution as $\tilde \XX$ started at $(x_1, \ldots, x_m)$.
\end{itemize}
We then set
\begin{equation*}
S_t = \{\phi(s,t,x) : (s,x) \in \Pi^+, \, 0 \le s \le t\},
\end{equation*}
where we use the short-hand notation $(s,x) \in \Pi^+$ to mean that
$(s,x)$ is an atom of $\Pi^+$.

For any $b>0$ we have that
almost surely
\begin{equation*}
S_t \cap [-b,b] =  \{\phi(s,t,x)\in [-b,b] : (s,x) \in \Pi^+, \, 0
\le s \le t, x \in [-a,a] \}
\end{equation*}
for all $a>0$ sufficiently large:  in particular, $S_t$ is almost
surely a discrete set and we identify $S_t$ interchangeably with
the simple point process obtained by placing a unit mass on each
point.  Using this observation, conditioning on $\Pi^+$, by
Proposition \ref{Pro2b}, and taking limits, we get the following
result which characterizes the {\em avoidance function} and hence
the distribution of $S_t$ (see Theorem 3.3 of \cite{Kal76}).

\begin{Prop}
\label{Pro:avoidance}
Given $y_1<\ldots<y_{2n}$, let $\YY$ be a coalescing Brownian
motion starting from $(y_1,\ldots,y_{2n})$. Then
\begin{equation*}
\begin{split}
\bP\left\{S_t\cap \bigcup_{j=1}^n \, ]y_{2j-1},y_{2j}] =
\emptyset\right\} =\bP\left[\exp\left(-\lambda\sum_{j=1}^n
\int_0^t Y_{2j}(s)- Y_{2j-1}(s)ds\right)\right].
\end{split}
\end{equation*}
\end{Prop}

We can re-phrase Proposition \ref{Pro:avoidance} as follows. For
fixed $a_0 \in \bR$, the function
\begin{equation*}
b \mapsto |D_t^\rightarrow(\phi(0,\cdot,b), \phi(0,\cdot,a_0))|
= \int_0^t \phi(0,s,b) - \phi(0,s,a_0) \, ds, \quad b \ge a_0,
\end{equation*}
where $|\cdot|$ denotes Lebesgue measure in the plane,
is non-negative, non-decreasing, and right-continuous.
It follows that there is a unique random Radon measure
$M_t$ on $\bR$ such that
\begin{equation*}
M_t(]a,b]) = \int_0^t \phi(0,s,b) - \phi(0,s,a) \, ds, \quad b \ge a.
\end{equation*}
Proposition \ref{Pro:avoidance} then says that $S_t$ is the simple
point process obtained by placing a unit mass at each atom of the
Cox process with the random intensity measure $\lambda M_t$; that
is, conditional on $M_t = m$, $S_t$ is distributed as the random
measure which places a unit mass at each atom of a Poisson process
with intensity measure $\lambda m$. Note that $M_t$ has atoms, and
so the resulting Cox process will not be a simple point process;
that is, it can have atoms with mass greater than one.
Consequently, $S_t$ is not a Cox process.

There is  a unique random Radon measure $M_\infty$ on $\bR$ such
that
\begin{equation*}
\begin{split}
M_\infty(]a,b])
& = \int_0^\infty \phi(0,s,b) - \phi(0,s,a) \, ds \\
& = \lim_{t \rightarrow \infty}  M_t(]a,b]) \\
& = \sup_{t \ge 0} M_t(]a,b]), \quad b \ge a, \\
\end{split}
\end{equation*}
(the finiteness of $M_\infty(]a,b])$ is assured by
the continuity of $s \mapsto \phi(0,s,a)$
and $s \mapsto \phi(0,s,b)$ and the fact that
$\phi(0,s,a) = \phi(0,s,b)$ for all $s$ sufficiently large).
Hence $S_t$ converges in distribution as $t \rightarrow \infty$
to the  simple point process obtained by placing a unit mass
at each atom of the
Cox process with the random intensity measure $\lambda M_\infty$.
As with $S_t$, the point process $S_\infty$ is not a Cox process.

We can give an almost sure construction of $S_\infty$
as follows. Consider a
Poisson random measure $\Pi^-$ on $\bR_- \times \bR$ with intensity
measure $\lambda \times {\mathrm{Lebesgue}}$.  Then
$S_t$ has the same distribution as
\begin{equation*}
\{\phi(s,0,x) : (s,x) \in \Pi^-, \, -t \le s \le 0\},
\end{equation*}
and so $S_\infty$ has the same distribution as
\begin{equation*}
\{\phi(s,0,x) : (s,x) \in \Pi^-, \, -\infty < s \le 0\}.
\end{equation*}

We can do some explicit computations for
$S_\infty$.  In what follows, let $\mathrm{Ai}$ denote
the Airy function -- see \cite{AbSt72} for its definition and
related properties.

\begin{Prop}
For $a<b$,
\[
\bP\{S_\infty\cap ]a,b]=\emptyset\}
=\frac{\mathrm{Ai}\left(\lambda^\frac{1}{3}(b-a)\right)}{\mathrm{Ai}(0)}
\]
and
\[
\bP[\# \, S_\infty\cap ]a,b]] = (3 \lambda)^{\frac{1}{3}}
\frac{\Gamma(\frac{2}{3})}{\Gamma(\frac{1}{3})} (b-a).
\]
\end{Prop}

\begin{proof}
Let $(Y_1,Y_2)$ be a two-dimensional coalescing Brownian motion
starting at $(a,b)$. Then $\frac{1}{\sqrt{2}}(Y_2-Y_1)$ is a
Brownian motion stopped at $0$. By Theorem 1 in \cite{Le89},
Theorem 1 in \cite{La93}, or Proposition 5.14 in \cite{GeSh}, we
have
\begin{equation*}
\begin{split}
\bP\{S_\infty\cap]a,b]=\emptyset\}
& = \bP\left[\exp(- \lambda M_\infty(]a,b]))\right] \\
&=\bP\left[\exp\left(-\lambda\int_0^\infty
 Y_2(s)-Y_1(s) \, ds \right)\right]\\
&=\frac{\mathrm{Ai}\left(\lambda^\frac{1}{3}(b-a)\right)}{\mathrm{Ai}(0)}.\\
\end{split}
\end{equation*}

Note that
\begin{equation*}
\begin{split}
\lim_{d - c \downarrow 0} \frac{\bP\{S_\infty \cap ]c,d]
\neq\emptyset\}}{d-c}
& =
- \frac{d}{dx}
\frac{\mathrm{Ai}\left(\lambda^\frac{1}{3} x \right)}{\mathrm{Ai}(0)}
\Bigg |_{x = 0} \\
& =
- \lambda^\frac{1}{3} \frac{\mathrm{Ai}'(0)}{\mathrm{Ai}(0)} \\
& =
(3 \lambda)^{\frac{1}{3}}
\frac{\Gamma(\frac{2}{3})}{\Gamma(\frac{1}{3})}. \\
\end{split}
\end{equation*}
Thus,
\begin{equation*}
\begin{split}
\bP[\# \, S_\infty\cap]a,b]]
&=\lim_{n\goto\infty}\sum_{i=1}^n
\bP\left\{S_\infty\cap\left]a+\frac{(i-1)(b-a)}{n},a+\frac{i(b-a)}{n}\right]\neq\emptyset\right\}\\
&=
(3 \lambda)^{\frac{1}{3}}
\frac{\Gamma(\frac{2}{3})}{\Gamma(\frac{1}{3})}
(b-a).\\
\end{split}
\end{equation*}
\end{proof}

\newcommand{\etalchar}[1]{$^{#1}$}

\end{document}